\documentclass[10pt]{amsart}
\usepackage{amsfonts,amssymb,amscd,amsmath,enumerate,verbatim,calc}
\usepackage{color,soul}
\usepackage{mathtools}
\usepackage{graphicx}
\everymath{\displaystyle}
\usepackage{fancyhdr}
\usepackage{tikz}
\usepackage{url}
\usetikzlibrary{shapes.geometric, arrows}
\tikzstyle{arrow} = [thick,->,>=stealth]
\tikzstyle{process} = [rectangle, minimum width=4cm, minimum height=2cm, text centered, text width=2.5cm, draw=green, fill=green!10]
\theoremstyle{plain}

\theoremstyle{definition}

\begin{document}
\title{SC*-Regular spaces and some functions}
\author{Neeraj kumar Tomar}
\email{neer8393@gmail.com}
\address{Department of Applied Mathematics, Gautam Buddha University, Greater Noida, Uttar Pradesh 201312, India}

\author{Amit Ujlayan}
\email{amitujlayan@gbu.ac.in}
\address{Department of Applied Mathematics, Gautam Buddha University, Greater Noida, Uttar Pradesh 201312, India}

\author{M. C. Sharma}
\email{sharmamc2@gmail.com}
\address{Prof. Department of Mathematics, N. R. E. C. College Khurja, Uttar Pradesh 203131, India}

\date{\today}
	
\begin{abstract}
This paper introduces a novel class of topological spaces, termed $SC^*$-regular spaces, which are defined using $SC^*$-open sets. We explore their fundamental properties and examine their connections with existing regularity concepts, such as regular, almost, softly, weakly, alpha, zeta, and generalized-regular spaces respectively. Furthermore, we  examined and define, analyze generalized $SC^*$-closed sets and $SC^*$-generalized closed functions, establishing key properties and preservation theorems. Several characterizations of $SC^*$-regular spaces are also presented, providing new insights into generalized regularity in topology.
\end{abstract}
\maketitle

\section*{Introduction} 

O. Njastad\cite{njȧstad1965some} introduced and examined the concept of $\alpha$-open sets, while M. K. Singal and S. P. Arya \cite{singal1969almost} defined two new classes of regular spaces: almost regular and weakly regular spaces. S. S. Benchalli \cite{benchalli2010some} introduced and studied the notion of $\alpha$-regular spaces. M. C. Sharma, P. Sharma and M. Singh \cite{sharma2014xi} introduced a new class of regular spaces called $\zeta$-regular spaces. H. Kumar \cite{kumar2017xi} obtained some more characterizations and preservation theorems for $\zeta$-regular spaces. H. Kumar and M. C. Sharma \cite{kumar2018softly} introduced two new classes of separation axioms, namely softly regular and partly regular spaces which are weaker than regular spaces. In 2015 Nidhi Sharma, et.al. \cite{Tomar2004normal} introduced the concept of $H^*$-open set and obtain some characterizations of $H^*$-normal spaces.\\
This paper is structured into five sections, each dedicated to different aspects of $SC^*$-regular spaces, their properties, and their relationships with other regularity concepts.\\
\textbf{Section 1.} This section introduces and defines generalized $SC^*$-closed and open sets, including $gSC^*$-closed, $SC^*g$-closed, $gSC^*$-open, and $SC^*g$-open sets. We present examples illustrating these concepts and explore their relationships with other known classes of open sets.\\
\textbf{Section 2.} In this section, we introduce $SC^*$-regular spaces using $SC^*$-open sets and examine their properties. We explore their relationships with regular, strongly $rg$-regular, almost regular, softly regular, weakly regular, and $\alpha$-regular spaces, along with their connections to $SC^*$-$T_1$, $SC^*$-$T_2$, and $SC^*$-compact Hausdorff spaces.\\
\textbf{Section 3.} This section examines functions related to $SC^*$-regular spaces. We introduce key theorems and lemmas concerning these functions and investigate their connections to strongly $SC^*$-open, strongly $SC^*$-closed, and almost $SC^*$-irresolute functions.  Additionally, we examine how these functions preserve $SC^*$-regularity. \\
\textbf{Section 4.} Here, we explore generalized $SC^*$-closed functions in the context of $SC^*$-regular spaces. We define and analyze relationships among various function types, including $gSC^*$-closed, $SC^*g$-closed, quasi $SC^*$-closed, $SC^*$-$SC^*g$-closed, $SC^*$-$gSC^*$-closed, and almost $gSC^*$-closed functions. Additionally, we present fundamental properties and a key theorem based on these functions.\\
\textbf{Section 5.} This section establishes further characterizations and preservation theorems for $SC^*$-regular spaces, presenting key results through theorems, lemmas, and corollaries. We also explore conditions under which $SC^*$-regularity is maintained under different types of mappings and investigate its behavior in product and subspace topologies.

\section{Preliminaries and Notations} 

In this paper, we consider topological spaces $(X, \tau)$, $(Y, \sigma)$, and $(Z, \gamma)$, where no separation axioms are assumed unless explicitly specified. Functions between these spaces are represented as $f: X \to Y$ and $g: Y \to Z$. For any subset $A$ of $X$, the notations $cl(A)$ and $int(A)$ denote its closure and interior, respectively.

We now define some basic notions which will be used throughout. For a good understanding of them, readers are referred to see \cite{levine1963semi, malathi2017pre,njȧstad1965some,stone1937applications}.

\subsection{Definition:} A subset $A$ of a topological space $X$ is said to be:
\begin{enumerate}
    \item  \textbf{regular open} \cite{stone1937applications} if $A=int(cl(A))$.
\item  \textbf{semi open} \cite{levine1963semi} if $A\subset cl(int(A))$.
\item  \textbf{$\alpha$-open} \cite{njȧstad1965some} if $A\subset (int(cl(int(A)))$.
\item  \textbf{$c^*$-open} \cite{malathi2017pre} if $int(cl(A))\subset A\subset cl(int(A))$.
\end{enumerate}

A set is termed \textbf{regular-closed} if its complement is regular-open. Similarly, the complements of semi-open, $\alpha$-open and $c^*$-open sets are referred to as \textbf{semi-closed, $\alpha$-closed}, and \textbf{$c^*$-closed}, respectively.

The \textbf{$c^*$-closure} and \textbf{semi-closure} of a set $A$, written as \textbf{$c^*$-$cl(A)$} and \textbf{$s$-$cl(A)$} respectively, are the smallest $c^*$-closed and semi-closed sets that contain $A$. Likewise, the \textbf{$c^*$-interior} and \textbf{semi-interior}, denoted \textbf{$c^*$-$int(A)$} and \textbf{$s$-$int(A)$}, represent the largest $c^*$-open and semi-open sets contained in $A$.

\subsection{Definition} A subset $A$ of a topological space $X$ is said to be \textbf{$SC^*$-closed}\cite{Chandrakala2022} if $scl(A)\subset U$ whenever $A\subset U$ and $U$ is $c^*$-open in $X$.\\
The complement of $SC^*$-closed set is said to be \textbf{$SC^*$-open}\cite{Chandrakala2022}.

The \( SC^* \)-closure of a set \( A \), denoted \( SC^*\text{-}cl(A) \), is defined as the smallest \( SC^* \)-closed set containing \( A \), obtained by intersection of all \( SC^* \)-closed supersets of \( A \). Similarly, \( SC^* \)-interior of \( A \), denoted \( SC^*\text{-}int(A) \), is the largest \( SC^* \)-open subset of \( A \), formed by the union of all \( SC^* \)-open sets contained in \( A \).

Collection of all $SC^*$-open, $SC^*$-closed, \(r\)-open, \(r\)-closed, \(s\)-open, and \(s\)-closed sets in a space $X$ is represented as $SC^*O(X)$, $SC^*C(X)$, $RO(X)$, $RC(X)$, $SO(X)$, and $SC(X)$, respectively.

\subsection{Definition} A subset $A$ of a topological space $(X,\tau )$ is said to be: 
\begin{enumerate}
    \item \textbf{$g$-closed}\cite{levine1970gen}  if $cl(A)\subset U$ whenever $A\subset U$ and 
$U\in \tau$. 
\item  \textbf{generalized $SC^*$-closed} (briefly $gSC^*$-closed)\cite{tomar2024scnormal} if  $SC^*$-$cl(A)\subset U$ whenever $A\subset U$ and $U$ is open in $X$.
\item  \textbf{$SC^*$generalized-closed} (briefly $SC^*g$-closed)\cite{tomar2024scnormal} if  $SC^*$-$cl(A)\subset U$ whenever $A\subset U$ and $U$ is $SC^*$-open.
\end{enumerate}

A set is termed \textbf{$g$-open, $gSC^*$-open}, or \textbf{$ SC^*g$-open} if its complement is $g$-closed, $gSC^*$-closed, or $SC^*g$-closed, respectively.
\vspace{.1cm}
\subsection{Remark.} The following relationships exist among different types of closed sets:

\vspace{8.0mm}

\begin{flushleft}
    closed $\Rightarrow$ $SC^*$-closed  $\Leftrightarrow$ $SC^*g$-closed  $\Leftrightarrow$  $gSC^*$- closed \\$\Downarrow$\\ $g$-closed
\end{flushleft}

\vspace{5.5mm}

However, these implications are generally one-way, as the reverse statements do not always hold. This can be illustrated with a suitable example.

\subsection{Example} Let $X=\{k,l,m,n\}$ and $\tau = \{\phi,\{k\},\{l\},\{k,l\},\{k,m\},\{k,l,m\},X\}$. Then

\begin{flushleft}
    $(1)$ closed sets in $(X,\tau)$ are $\phi,\{n\},\{l,n\},\{m,n\},\{k,m,n\},\{l,m,n\},X$

$(2)$ $g$-closed set in $(X,\tau)$ are $\phi,X,\{n\},\{k,n\},\{l,n\},\{m,n\},\{k,l,n\},\{k,m,n\},\{l,m,n\}$.

$(3)$ $SC^*$-closed set in $(X,\tau)$ are $\phi$, $X$, $\{k\}$, $\{l\}$, $\{m\}$, $\{n\}$, $\{k,l\}$, $\{k,m\}$, $\{k,n\}$, $\{l,m\}$, $\{l,n\}$, $\{m,n\}$, $\{k,l,m\}$, $\{k,l,n\}$, $\{k,m,n\}$, $\{l,m,n\}$.

$(4)$ $gSC^*$-closed set in $(X,\tau)$ are $\phi$, $X$, $\{k\}$, $\{l\}$, $\{m\}$, $\{n\}$, $\{k,l\}$, $\{k,m\}$, $\{k,n\}$, $\{l,m\}$, $\{l,n\}$, $\{m,n\}$, $\{k,l,m\}$, $\{k,l,n\}$, $\{k,m,n\}$, $\{l,m,n\}$.

$(5)$ $SC^*g$-closed set in $(X,\tau)$ are $\phi$, $X$, $\{k\}$, $\{l\}$, $\{m\}$, $\{n\}$, $\{k,l\}$, $\{k,m\}$, $\{k,n\}$, $\{l,m\}$, $\{l,n\}$, $\{m,n\}$, $\{k,l,m\}$, $\{k,l,n\}$, $\{k,m,n\}$, $\{l,m,n\}$.
\end{flushleft}

\subsection{Lemma.} 
If \( J \subset X\) where, \(X\) is a topological space and let \( SC^*\text{-}cl(J) \) denote its \( SC^* \)-closure. The following statements hold:\\
    (i) A point \( x \in X \) belongs to \( SC^*\text{-}cl(J) \) iff every \( SC^* \)-open set containing \( x \) intersects \( J \), i.e., \( J \cap U \neq \emptyset \) for all \( U \in SC^*O(X) \) with \( x \in U \).\\
    (ii) The set \( J \) is \( SC^* \)-closed precisely when \( J = SC^*\text{-}cl(J) \).\\
    (iii) For any subsets \( J \subseteq I \subseteq X \), it follows that \( SC^*\text{-}cl(J) \subseteq SC^*\text{-}cl(I) \).\\
    (iv) The \( SC^* \)-closure operator is idempotent, meaning \( SC^*\text{-}cl(SC^*\text{-}cl(J)) = SC^*\text{-}cl(J) \).\\
    (v) The \( SC^* \)-closure of any subset is itself an \( SC^* \)-closed set.

\subsection{Lemma.} 
If \( J \subset X\) where, \(X\) is a topological space is said to be \( gSC^* \)-open iff, for every closed set \( F \subseteq J \), it follows that \( F \subseteq SC^*\text{-}int(J) \).

 \section{$SC^*$-regular spaces}
\subsection{Definition} A space $X$ is said to be regular (resp. $g$-regular\cite{munshi1986separation}, $SC^*$-regular) if for every closed set $F$ and a point $x$ $\notin$ $F$, there exist disjoint open (resp. $g$-open, $SC^*$-open) sets $U$ and $V$ of $X$ such that $F\subset U$ and $x\in V$.
Clearly every regular space is $SC^*$-regular space.

\subsection{Definition} A space $X$ is said to be softly regular \cite{kumar2018softly} (resp. almost regular \cite{singal1969almost}, strongly 
$rg$-regular \cite{gnanachandra2011strongly}) if for every $\pi$-closed (resp. regular closed, $rg$-closed) set $F$ of $X$, and a point $x \in X- F$, there exist disjoint open sets $U$ and $V$ such that $F \subset U$ and $x \in V$.

\subsection{Definition} A space $X$ is said to be weakly regular \cite{singal1969almost} if for every point $x$ and every regularly open set $U$ containing $x$, there is an open set $V$ such that $x \in V \subset cl(V) \subset U$.

\subsection{Remark.} Hierarchy of Regularity in a Topological Space $(X,\tau)$:
\vspace{3.7mm}

\begin{flushleft}
    strongly $rg$-regular $\Rightarrow$ regular $\Rightarrow$ $\alpha$ -regular  $\Rightarrow$  $SC^*$-regular \\$\Downarrow$\\ softly regular $\Rightarrow$ almost regular $\Rightarrow$  weakly regular
\end{flushleft}

\vspace{3.9mm}

This diagram illustrates the logical progression of different regularity concepts in topological spaces, showing how one property implies another.

\subsection{Example} Let $X$ = $\{k,l,m,n\}$ and $\tau$ = $\{\phi,\{k\},\{l\},\{k,l\},\{m,n\},\{k,m,n\},\{l,m,n\},X\}$. Then $A$ = $\{k\}$ and $B$ = $\{\phi\}$ are disjoint closed sets,  $U$ = $\{k,m,n\}$ and  $V$ = $\{l\}$ such that $A\subset U$ and $B\subset V$. Hence $X$ regular, $g$-regular and as well as $SC^*$-regular because every regular is $SC^*$-regular.

\subsection{Example} If the topological space $X = \{k,l,m\}$ with the topology $\tau = \{\emptyset, \{k\}, \{l\}, \{k,l\}, X\}$. In this case, $X$ satisfies the conditions for weak regularity. However, it does not fulfill the criteria for almost regularity or soft regularity.

\subsection{Example} If $X = \{l, m, c, d\}$ with topology $\tau = \{\emptyset, \{k\}, \{l\}, \{k,l\}, \{k,l,m\}, \{k,l,n\}, X\}$. This space satisfies the conditions for almost regularity. However, it does not meet the requirements for strong $rg$-regularity.

\subsection{Example} If $X$ = $\{k,l,m\}$ with topology $\tau$ = $\{\phi,\{k\},\{l\},\{k,l\},\{k,m\},X\}$. Then the space $X$ is regular.

\subsection{Example} If $X = \{k, l, m\}$ with the topology $\tau = \{\emptyset, \{k\}, \{l,m\}, X\}$. This space satisfies the conditions for regularity. However, it does not fulfill the requirements for strong $rg$-regularity. Specifically, the set $F = \{l\}$ is $rg$-closed, but the point $m \notin F$. It is not possible to find disjoint open sets that separate $\{m\}$ and $\{l\}$.

\subsection{Theorem} The following statements are equivalent for a topological space $X$:
\begin{flushleft}
    $(i)$ The space $X$ satisfies the $SC^*$-regularity condition.\\
    $(ii)$ For every point \(x\) in \(  X \) and every open set \( M \) containing \( x \), there exists a set \( N \in SC^*O(X) \) such that:
    \[
    x \in N \subseteq SC^*\text{-cl}(N) \subseteq M.
    \]\\
    $(iii)$ Every closed subset \( F \subseteq X \) satisfies:
    \[
    F = \bigcap \{ SC^*\text{-cl}(N) \mid F \subseteq N,\ N \in SC^*O(X) \}.
    \]\\
    $(iv)$ Given any subset \( J \subseteq X \) and any open set \( M \) with \( J \cap M \neq \emptyset \), there exists \( N \in SC^*O(X) \) such that:
    \[
    J \cap N \neq \emptyset \quad \text{and} \quad SC^*\text{-cl}(N) \subseteq M.
    \]\\
    $(v)$ For any nonempty set \( J \subseteq X \) and any closed set \( F \subseteq X \) with \( J \cap F = \emptyset \), there exist \( N, W \in SC^*O(X) \) such that:
    \[
    J \cap N \neq \emptyset,\quad F \subseteq W,\quad \text{and} \quad N \cap W = \emptyset.
    \]

\end{flushleft}  
\vspace{1.7mm}

    \textbf{Proof:}  
\textbf{$(i)$ $\Rightarrow$ $(ii)$:}  
Suppose if $M$ be an open set containing $x$. Since $X - M$ is closed, does not include $x$, we use condition \textbf{$(i)$}. Thus, $\exists$ sets $W, N \in SC^*O(X)$ s.t.:  
\[
    x \in N, \quad X - M \subset W, \quad N \cap W = \emptyset.
\]  
By using \textbf{Lemma 1.9}, it follows that:  
\[
    SC^*\text{-cl}(N) \cap W = \emptyset.
\]  
Therefore,  
\[
    x \in N \subset SC^*\text{-cl}(N) \subset M.
\]  

\textbf{$(ii)$ $\Rightarrow$ $(iii)$:}  
Suppose that $F$ is closed subset of $X$. Then $F \subset N$, Thus using \textbf{Lemma 1.9 $(iii)$}, we obtain:  
\[
    SC^*\text{-cl}(F) \subset SC^*\text{-cl}(N),
\]  
Thus implies:  
\[
    F \subset SC^*\text{-cl}(N),
\]  
since $F \subset SC^*\text{-cl}(F)$. Therefore,  
\[
    \bigcap \{ SC^*\text{-cl}(N) : F \subset N \in SC^*O(X) \} \supseteq F.
\]  

For the reverse inclusion, assume $x \notin F$. So $X - F$ is open set and containing $x$. By \textbf{$(ii)$}, there exists $M \in SC^*O(X)$ s.t.:  
\[
    x \in M \subset SC^*\text{-cl}(M) \subset X - F.
\]  
Define $N = X - SC^*\text{-cl}(M)$. Again by using \textbf{Lemma 1.9}, we get:  
\[
    F \subset N \in SC^*O(X), \quad x \notin SC^*\text{-cl}(N).
\]  
Implies that:  
\[
    \bigcap \{ SC^*\text{-cl}(N) : F \subset N \in SC^*O(X) \} \subseteq F.
\]  
Thus,  
\[
    \bigcap \{ SC^*\text{-cl}(N) : F \subset N \in SC^*O(X) \} = F.
\]  

\textbf{$(iii)$ $\Rightarrow$ $(iv)$:}  
Let $J\subset X$ and suppose $M$ is an open set with $J \cap M \neq \emptyset$.  
For any $x \in J \cap M$, the set $X - M$ is closed and does not contain $x$.  

Using \textbf{$(iii)$}, $\exists$ $W \in SC^*O(X)$ s.t.:  
\[
    X - M \subset W, \quad x \notin SC^*\text{-cl}(W).
\]  
Define $N = X - SC^*\text{-cl}(W)$. So:  
\[
    N \subset X - W.
\]  
Since $x \in N \cap J$, applying \textbf{Lemma 1.9}, we obtain:  
\[
    N \in SC^*O(X), \quad SC^*\text{-cl}(N) \subset SC^*\text{-cl}(X - W) = X - W \subset M.
\]  

\textbf{$(iv)$ $\Rightarrow$ $(v)$:}  
Let $J\subset X$ and $F$ is closed set s.t. $J \cap F = \emptyset$. Since $X - F$ is an open set containing $J$, we apply condition \textbf{$(iv)$}.  
There exists $N \in SC^*O(X)$ such that:  
\[
    J \cap N \neq \emptyset, \quad SC^*\text{-cl}(N) \subset X - F.
\]  
Define $W = X - SC^*\text{-cl}(N)$. Then:  
\[
    F \subset W, \quad N \cap W = \emptyset.
\]  
By \textbf{Lemma 1.9}, we also conclude that $W \in SC^*O(X)$.  

\textbf{$(v)$ $\Rightarrow$ $(i)$:}  
This result follows directly from the given conditions, completing the proof. 

\subsection{Theorem}
For a topological space \(X\), \(X\) is \(SC^*\)-regular iff for every closed set \(F\in X\) and each point \(x \in X \setminus F\), $\exists$ \(SC^*\)-open sets \(M\) and \(N\) s.t. \(x \in M\), \(F \subset N\), and  
\[
SC^*\text{-cl}(M) \cap SC^*\text{-cl}(N) = \emptyset.
\]

\textbf{Proof:}    
\textbf{(Necessity)} Assume that \(X\) is \(SC^*\)-regular. Suppose that \(F\) be a closed subset of \(X\) and let \(x \notin F\). By the definition of \(SC^*\)-regularity, $\exists$ \(SC^*\)-open sets \(M_x\), \(N\) s.t. \(x \in M_x\), \(F \subset N\), and \(M_x \cap N = \emptyset\). So implies  
\[
M_x \cap SC^*\text{-cl}(N) = \emptyset.
\]  
If \(SC^*\text{-cl}(N)\) is \(SC^*\)-closed and does not contain \(x\), there must exist \(SC^*\)-open sets \(G_1\) and \(G_2\) such that \(x \in G_1\), \(SC^*\text{-cl}(N) \subset G_2\), and \(G_1 \cap G_2 = \emptyset\). Consequently, we get  
\[
SC^*\text{-cl}(G_1) \cap G_2 = \emptyset.
\]  
Now, let \(M = M_x \cap G_1\). Since both \(M_x\) and \(G_1\) are \(SC^*\)-open, their intersection is also \(SC^*\)-open. Thus, \(M\) and \(N\) satisfy \(x \in M\), \(F \subset N\), and  
\[
SC^*\text{-cl}(M) \cap SC^*\text{-cl}(N) \subseteq SC^*\text{-cl}(G_1) \cap G_2 = \emptyset.
\]  

\textbf{(Sufficiency)} Conversely, assume that for any closed set \(F\in X\) and any point \(x \notin F\), $\exists$ \(SC^*\)-open sets \(M\), \(N\) s.t. \(x \in M\), \(F \subset N\), and  
\[
SC^*\text{-cl}(M) \cap SC^*\text{-cl}(N) = \emptyset.
\]  
Since  
\[
M \cap N \subset SC^*\text{-cl}(M) \cap SC^*\text{-cl}(N),
\]  
we conclude that  
\[
M \cap N = \emptyset.
\]  
Thus, the given condition ensures that \(X\) satisfies the definition of \(SC^*\)-regularity, completing the proof.

\subsection{Definition}
A topological space \( X \) is called an \( SC^* \)-\( T_3 \) space if it satisfies both the \( SC^* \)-regularity condition and the \( SC^* \)-\( T_1 \)\cite{tomarsome2024} separation axiom.

\subsection{Theorem}
Every \( SC^* \)-\( T_3 \) space is necessarily an \( SC^* \)-\( T_2 \) space.

\textbf{Proof}. Assume that \(X\) is an \(SC^*\)-\(T_3\) space, meaning it satisfies both the \(SC^*\)-\(T_1\) and \(SC^*\)-regularity conditions.  

Since \(X\) is \(SC^*\)-\(T_1\), each singleton set \(\{x\}\) is an \(SC^*\)-closed set. Let \(x, y \in X\) be distinct points. Because \(X\) is \(SC^*\)-regular, $\exists$ two disjoint \(SC^*\)-open sets \(M\), \(N\) respectively so,  
\[
\{x\} \subset M, \quad y \in N, \quad M \cap N = \emptyset.
\]
This shows that \(x\) and \(y\) have disjoint \(SC^*\)-open neighborhoods.  

Since this holds for any two distinct points in \(X\), the space satisfies the \(SC^*\)-\(T_2\) (Hausdorff) separation property.  

Therefore, \(X\) is an \(SC^*\)-\(T_2\) space.

\subsection{Theorem}
Any subspace of an \( SC^* \)-regular topological space is itself \( SC^* \)-regular.

\textbf{Proof}: If \( X \) be an \( SC^* \)-regular space and let \( Y \subseteq X \) be a subspace. Suppose \( x \in Y \), \( F \) is a closed subset of \( Y \) s.t. \( x \notin F \).  

Since the subspace topology on \( Y \) is inherited from \( X \), there exists a closed set \( J \subseteq X \) such that \( F = Y \cap J \), and clearly \( x \notin J \).  

Because \( X \) is \( SC^* \)-regular, $\exists$ sets \( M, N \in SC^*O(X) \) s.t.:
\[
x \in M, \quad J \subseteq N, \quad \text{and} \quad M \cap N = \emptyset.
\]

Now, define the sets:
\[
M' = M \cap Y, \quad N' = N \cap Y.
\]
Then \( M', N' \in SC^*O(Y) \), since the intersection of an \( SC^* \)-open set in \( X \) with \( Y \) is \( SC^* \)-open in the subspace topology.

It follows that:
\[
x \in M', \quad F \subseteq N', \quad M' \cap N' = \emptyset.
\]

Hence, \( Y \) satisfies the \( SC^* \)-regularity condition, and so \( Y \) is \( SC^* \)-regular.

\subsection{Theorem}  
Each \( SC^* \)-comp., Hausdorff space is an space \( SC^* \)-\( T_3 \), consequently, an \( SC^* \)-regular space.

\textbf{Proof}. Suppose a comp., Hausdorff space \(X\), which implies it's a \(SC^*\)-\(T_3\) space. Since every \(SC^*\)-\(T_2\) space is also \(SC^*\)-\(T_1\), we need to establish that \(X\) satisfies the \(SC^*\)-regularity condition.  

Consider a closed subset \(F\) of \(X\) and a point \(x\) that is not in \(F\). If \(X\) is an Hausdorff space, for each point \(y\) in \(F\), $\exists$ two disjoint \(SC^*\)-open sets \(J_y\), \(G_y\) s.t. \(x \in J_y\) and \(y \in G_y\).  

Since \(F\) is a closed subset of \(X\), we consider the induced topology on \(F\), denoted as \(\zeta^*\). The collection  
\[
C^* = \{F \cap J_y : y \in F\}
\]
forms a covering of \(F\) by \(SC^*\)-open sets. Since \(X\) is \(SC^*\)-compact, the subspace \(F\) is also \(SC^*\)-compact. Therefore, a finite subcover exists, meaning there are points \(y_1, y_2, \dots, y_n \in F\) such that  
\[
F = \bigcup_{i=1}^{n} (F \cap J_{y_i}).
\]

Define \(J\) as  
\[
J = \bigcup_{i=1}^{n} J_{y_i}.
\]
Since each \(J_{y_i}\) is \(SC^*\)-open, their union \(J\) is also \(SC^*\)-open and contains \(F\).  

Similarly, since each \(G_{y_i}\) is \(SC^*\)-open and contains \(x\), their intersection  
\[
G = \bigcap_{i=1}^{n} G_{y_i}
\]
is also \(SC^*\)-open and contains \(x\).  

Moreover, \(G\) and \(J\) are disjoint, as otherwise, some \(G_{y_i}\) and \(J_{y_i}\) would intersect, contradicting their original disjointness.  

Thus, $\forall$ closed set \(F\) and point \(x\) outside \(F\), then we have disjoint \(SC^*\)-open sets \(G\), \(J\) respectively s.t. \(x \in G\) and \(F \subseteq J\). This verifies that \(X\) is \(SC^*\)-regular. Since \(X\) is also \(SC^*\)-\(T_2\), it follows that \(X\) is an \(SC^*\)-\(T_3\) space.

\section{Function Classes Related to \( SC^* \)-Regular Spaces}

\subsection{Definition} A function $f : X \rightarrow Y $ is called.

\begin{flushleft}
    $(1)$ \textbf{R-map}\cite{carnahan1974some} if $f^{-1}(V)$ is regular open in $X$ for every regular open set $V$ of $Y$,
\end{flushleft}  

\begin{flushleft}
    $(2)$ \textbf{completely continuous}\cite{arya1974strongly} if $f^{-1}(V)$ is regular open in $X$ for every open set $V$ of $Y$,
\end{flushleft}

\begin{flushleft}
    $(3)$ \textbf{rc-continuous}\cite{jankovic1985note} if for each regular closed set $F$ in $Y, f^{-1}(F)$ is regular closed in $X$.
\end{flushleft}

\subsection{Definition} A mapping $X$ to $Y$ is described as follows:

\begin{flushleft}
    $(1)$ \textbf{strongly-$SC^*$-open}: If for every set $M$ in $SC^*O(X)$, its image $f(M)$ belongs to $SC^*O(Y)$.  
\end{flushleft}

\begin{flushleft}
    $(2)$ \textbf{strongly-$SC^*$-closed}: If for every set $M$ in $SC^*C(X)$, the image $f(M)$ is a member of $SC^*C(Y)$.  
\end{flushleft}

\begin{flushleft}
    $(3)$ \textbf{almost-$SC^*$-irresolute}: If any point $x$ in $X$, any $SC^*$-neighborhood $N$ of $f(x)$, the $SC^*$-closure of $f^{-1}(N)$, written as $SC^*$-$cl(f^{-1}(N))$, remains an $SC^*$-neighborhood of $x$.  
\end{flushleft}

\subsection{Theorem} A mapping \(X \) to \(Y\) is considered strongly-\( SC^* \)-closed iff, for every subset \( J \) of \( Y \) and any \( SC^* \)-open set \( M \) in \( X \) that contains \( f^{-1}(J) \), $\exists$ an \( SC^* \)-open set \( N \) in \( Y \) s.t. \( J \subseteq N \) and \( f^{-1}(N) \subseteq M \).  

\textbf{Proof}.  
\textbf{(Sufficiency)}:  
Suppose \( f \) is strongly-\( SC^* \)-closed. Given any subset \( J \) of \( Y \), let \( M \) be an \( SC^* \)-open set in \( X \) containing \( f^{-1}(J) \). Define  
\[
N = Y \setminus f(X \setminus M).
\]
Since the complement of an \( SC^* \)-closed set is \( SC^* \)-open, \( N \) is an \( SC^* \)-open set in \( Y \) that satisfies \( J \subseteq N \) and \( f^{-1}(N) \subseteq M \), as required.  

\textbf{(Necessity)}:  
Conversely, suppose that any subset \( J \) of \( Y \), any \( SC^* \)-open set \( M \) in \( X \) containing \( f^{-1}(J) \), $\exists$ an \( SC^* \)-open set \( N \) in \( Y \) s.t., \( J \subseteq N \) and \( f^{-1}(N) \subseteq M \).  

To demonstrate that \( f \) is strongly-\( SC^* \)-closed, consider any \( SC^* \)-closed set \( L \) in \( X \). Since \( L \) is \( SC^* \)-closed, its complement \( X \setminus L \) is \( SC^* \)-open. Also,  
\[
f^{-1}(Y \setminus f(L)) \subseteq X \setminus L.
\]
By the given condition, $\exists$ an \( SC^* \)-open set \( N \) in \( Y \) s.t.:  
\[
Y \setminus f(L) \subseteq N \quad \text{and} \quad f^{-1}(N) \subseteq X \setminus L.
\]
Rearranging these inclusions, we obtain  
\[
f(L) \supseteq Y \setminus N \quad \text{and} \quad L \subseteq f^{-1}(Y \setminus N).
\]
This leads to  
\[
f(L) = Y \setminus N.
\]
Since \( Y \setminus N \) is \( SC^* \)-closed, we conclude that \( f(L) \) is also \( SC^* \)-closed, proving that \( f \) is strongly-\( SC^* \)-closed.

\subsection{Lemma.} Let \(X \) to \(Y\) be a mapping. Then conditions are equivalent: \\
$(i)$ \( f \) is almost-\( SC^* \)-irresolute.\\
$(ii)$ $\forall$, \( N \) in \( SC^*O(Y) \), we have  
    \[
    f^{-1}(N) \subseteq SC^* \text{-} \operatorname{int} (SC^* \text{-} \operatorname{cl} (f^{-1}(N))).
    \]

\subsection{Theorem} A function \(X \) to \(Y\) is almost-\( SC^* \)-irresolute iff for every \( SC^* \)-open set \( M \) in \( X \), the inclusion  
\[
SC^* \text{-} cl(M) \subseteq SC^* \text{-} cl(f(M))
\]  
holds.

\textbf{Proof}.  
\textbf{(Necessity)}:  
Assume \( f \) is almost \( SC^* \)-irresolute, and let \( M \) be any \( SC^* \)-open set in \( X \). Suppose $\exists$ \( y \notin SC^* \text{-} cl(f(M)) \). Then,  $\exists$ an \( SC^* \)-open set \( N \) in \( Y \) s.t.:  
\[
N \cap f(M) = \emptyset.
\]
This implies that  
\[
f^{-1}(N) \cap M = \emptyset.
\]
Since \( M \) is \( SC^* \)-open, then  
\[
SC^* \text{-} int(SC^* \text{-} cl(f^{-1}(N))) \cap SC^* \text{-} cl(M) = \emptyset.
\]
using \textbf{Lemma 3.4}, it follows that  
\[
f^{-1}(N) \cap SC^* \text{-} cl(M) = \emptyset.
\]
Thus,  
\[
N \cap f(SC^* \text{-} cl(M)) = \emptyset,
\]
which leads to \( y \notin f(SC^* \text{-} cl(M)) \). Hence,  
\[
SC^* \text{-} cl(M) \subseteq SC^* \text{-} cl(f(M)).
\]

\textbf{(Sufficiency)}:  
Now, assume that for every \( SC^* \)-open set \( M \) in \( X \),  
\[
SC^* \text{-} cl(M) \subseteq SC^* \text{-} cl(f(M)).
\]
Let \( N \) be an \( SC^* \)-open set in \( Y \). Define  
\[
G = X \setminus SC^* \text{-} cl(f^{-1}(N)).
\]
Since \( SC^* \text{-} cl(f^{-1}(N)) \) is \( SC^* \)-closed, it follows that \( G \) is \( SC^* \)-open in \( X \). By the hypothesis,  
\[
f(SC^* \text{-} cl(G)) \subseteq SC^* \text{-} cl(f(G)).
\]
Rewriting this, we obtain  
\[
X \setminus SC^* \text{-} int(SC^* \text{-} cl(f^{-1}(N))) = SC^* \text{-} cl(G).
\]
Using the assumption,  
\[
SC^* \text{-} cl(G) \subseteq f^{-1}(SC^* \text{-} cl(f(G))) \subseteq f^{-1}(SC^* \text{-} cl(f(X \setminus f^{-1}(N)))).
\]
Since  
\[
f(X \setminus f^{-1}(N)) = Y \setminus N,
\]
we conclude  
\[
f^{-1}(SC^* \text{-} cl(Y \setminus N)) = f^{-1}(Y \setminus N) = X \setminus f^{-1}(N).
\]
Thus,  
\[
f^{-1}(N) \subseteq SC^* \text{-} int(SC^* \text{-} cl(f^{-1}(N))).
\]
By again using \textbf{Lemma 3.4}, this confirms that \( f \) is almost-\( SC^* \)-irresolute. Hence, the proof is complete.

\subsection{Theorem.} Suppose \( f: X \to Y \) is a surjective function that is strongly-\( SC^* \)-open, continuous, and almost-\( SC^* \)-irresolute. If \( X \) is an \( SC^* \)-normal space, then \( Y \) must be an \( SC^* \)-regular space.

\textbf{Proof}:  
Consider a closed set \( J \) in \( Y \) and an open set \( I \) in \( Y \) s.t. \( J \subset I \). Since \( f \) is cont., its preimage \( f^{-1}(J) \) closed in \( X \), \( f^{-1}(I) \) open in \( X \), which ensures, \( f^{-1}(J) \subset f^{-1}(I) \).  

If \( X \) is \( SC^* \)-normal, $\exists$ an \( SC^* \)-open set \( M \) in \( X \) s.t.:  
\[
f^{-1}(J) \subset M \subset SC^* \text{-} \operatorname{cl}(M) \subset f^{-1}(I).
\]  
Applying the properties of \( f \), we get  
\[
f(f^{-1}(J)) \subset f(M) \subset f(SC^* \text{-} \operatorname{cl}(M)) \subset f(f^{-1}(I)).
\]  
If \( f \) is strongly-\( SC^* \)-open and almost-\( SC^* \)-irresolute, it follows that  
\[
J \subset f(M) \subset SC^* \text{-} \operatorname{cl}(f(M)) \subset I.
\]  
Thus, the space \( Y \) satisfies the conditions for \( SC^* \)-regularity, completing the proof.

\subsection{Theorem.} Let \( f: X \to Y \) be a strongly-\( SC^* \)-closed and continuous function, where \( X \) is an \( SC^* \)-regular. Then \( Y \) is also \( SC^* \)-regular.

\textbf{Proof}: Suppose that \( G_1 \), \( G_2 \) two disjoint closed subsets of \( Y \). And \( f \) is cont., their preimages \( f^{-1}(G_1) \) and \( f^{-1}(G_2) \) are also closed in \( X \). Given that \( X \) is \( SC^* \)-regular, $\exists$ disjoint \( SC^* \)-open sets \( M \), \( N \) in \( X \) s.t. \( f^{-1}(G_1) \subset M \), \( f^{-1}(G_2) \subset N \). By \textbf{Theorem 3.3}, $\exists$ \( SC^* \)-open sets \( J \) and \( I \) in \( Y \) s.t. \( G_1 \subset J \), \( G_2 \subset I \), and their preimages satisfy \( f^{-1}(J) \subset M \) and \( f^{-1}(I) \subset N \). Since \( M \) and \( N \) are disjoint, it follows that \( J \) and \( I \) are also disjoint. This confirms that \( Y \) satisfies the conditions of \( SC^* \)-regularity.

\section{Generalized $SC^*$-closed functions and some theorems}

\subsection{Definition.} A function \( f: X \to Y \) is defined as follows:\\  
    $(1)$ \textbf{\( SC^* \)-closed}:\cite{Chandrakala2022} The image \( f(J) \) is \( SC^* \)-closed in \( Y \) for every closed subset \( J \) of \( X \).\\  
     $(2)$ \textbf{\( SC^*g \)-closed}: The image \( f(J) \) is \( SC^*g \)-closed in \( Y \) for every closed subset \( J \) of \( X \).\\  
     $(3)$ \textbf{\( gSC^* \)-closed}: The image \( f(J) \) is \( gSC^* \)-closed in \( Y \) for every closed subset \( J \) of \( X \).

\subsection{Definition.} A function \( f: X \to Y \) is classified as follows:\\  
     $(1)$ \textbf{quasi \( SC^* \)-closed}: The image \( f(J) \) is closed in \( Y \) for every \( J \in SC^*C(X) \).\\  
     $(2)$ \textbf{\( SC^* \)-\( SC^*g \)-closed}: The image \( f(J) \) is \( SC^*g \)-closed in \( Y \) for every \( J \in SC^*C(X) \).\\  
     $(3)$ \textbf{\( SC^* \)-\( gSC^* \)-closed}: The image \( f(J) \) is \( gSC^* \)-closed in \( Y \) for every \( J \in SC^*C(X) \).\\  
     $(4)$ \textbf{almost \( gSC^* \)-closed}:\cite{tomarsome2024} The image \( f(J) \) is \( gSC^* \)-closed in \( Y \) for every \( J \in RC(X) \).

\subsection{Definition.} A function \( f: X \to Y \) is called \textbf{\( SC^* \)-\( gSC^* \)-continuous} if for every\\ \( L \in SC^*C(Y) \), the preimage \( f^{-1}(L) \) is \( gSC^* \)-closed in function \( X \).

\subsection{Remark}  
Every closed function is inherently an \( SC^* \)-closed function; however, the reverse does not necessarily hold. Furthermore, since every \( SC^* \)-closed set is also \( gSC^* \)-closed, it follows that every \( SC^* \)-closed function is \( gSC^* \)-closed. Thus, both \( SC^* \)-closed and \( SC^* \)-\( gSC^* \)-closed functions satisfy the condition for being \( gSC^* \)-closed.

\subsection{Theorem}  
A surjective function \(X \) to \(Y\) is \( gSC^* \)-closed (respectively, \( SC^* \)-\( gSC^* \)-closed) iff for every subset \( I \subset Y \) and for every open (respectively, \( SC^* \)-open) set \( M \subset X \) with \( f^{-1}(I) \subset M \), $\exists$ a \( gSC^* \)-open set \( N \subset Y \) s.t. \( I \subset N \) and \( f^{-1}(N) \subset M \).

\textbf{Proof:}  
Assume \( f \) is \( gSC^* \)-closed (or \( SC^* \)-\( gSC^* \)-closed). Let \( I \subset Y \) and let \( M \) be an open (or \( SC^* \)-open) subset of \( X \) s.t. \( f^{-1}(I) \subset M \). Consider the set \( N = Y \setminus f(X \setminus M) \). Then,  
\[
N^c = Y \setminus N = f(X \setminus M).
\]  
Since \( X \setminus M \) is closed (or \( SC^* \)-closed) and \( f \) is \( gSC^* \)-closed (or \( SC^* \)-\( gSC^* \)-closed), the image \( f(X \setminus M) \) is \( gSC^* \)-closed in \( Y \). Hence, \( N \) is \( gSC^* \)-open, and clearly \( I \subset N \) with \( f^{-1}(N) \subset M \).

Conversely, suppose the stated condition holds. Let \( F \subset X \) be closed (or \( SC^* \)-closed), and define \( I = Y \setminus f(F) \). Then \( f^{-1}(I) \subset X \setminus F \), which is open (or \( SC^* \)-open). By hypothesis, $\exists$ a \( gSC^* \)-open set \( N \subset Y \) s.t. \( I \subset N \), \( f^{-1}(N) \subset X \setminus F \). It follows that \( N \subset Y \setminus f(F) = I \), so \( N = I \), and thus \( f(F) = Y \setminus N \), which is \( gSC^* \)-closed. Therefore, \( f \) is \( gSC^* \)-closed (or \( SC^* \)-\( gSC^* \)-closed), as required.

\subsection{Remark}  
Necessity part of preceding theorem can be derived by considering each set in its closed form, as elaborated in the following proposition.

\subsection{Proposition} Let \( f: X \to Y \) be a surjective function that is \( gSC^* \)-closed (or \( SC^* \)-\( gSC^* \)-closed). If \( F \) closed subset of \( Y \) $\&$ \( M \) an open (or \( SC^* \)-open) subset of \( X \) containing \( f^{-1}(F) \), $\exists$ an \( SC^* \)-open set \( N \) in \( Y \) s.t. \( F \subset N \), \( f^{-1}(N) \subset M \).  

\textbf{Proof:} According to \textbf{Theorem 4.5}, $\exists$ a \( gSC^* \)-open set \( W \subseteq Y \) s.t. \( F \subseteq W \) $\&$ \( f^{-1}(W) \subseteq M \). If \( F \) is closed in \( Y \), it follows from \textbf{Lemma 1.9} that \( F \subseteq SC^*\text{-}int(W) \), where \( SC^*\text{-}int(W) \) denotes the \( SC^* \)-interior of \( W \). Define \( N = SC^*\text{-}int(W) \). Then \( N \) is an \( SC^* \)-open subset of \( Y \), satisfying \( F \subseteq N \) and \( f^{-1}(N) \subseteq f^{-1}(W) \subseteq M \), as required.

\subsection{Proposition} 
Let a mapping \(X \) to \(Y\) be a cont., $\&$ \( SC^* \)-\( gSC^* \)-closed function. If \( J \subseteq X \) is \( gSC^* \)-closed, thus \( f(J) \) is \( gSC^* \)-closed in \( Y \).

\textbf{Proof:}  
Let \( N \subseteq Y \) be an open set s.t. \( f(J) \subseteq N \). Since \( f \) is cont., the preimage \( f^{-1}(N) \) is open in \( X \), clearly \( J \subseteq f^{-1}(N) \). Given \( J \) is \( gSC^* \)-closed, it follows by definition that
\[
SC^*\text{-}\mathrm{cl}(J) \subseteq f^{-1}(N).
\]
Applying \( f \) to both sides, we obtain
\[
f(SC^*\text{-}\mathrm{cl}(J)) \subseteq N.
\]
Since \( SC^*\text{-}\mathrm{cl}(J) \) is \( SC^* \)-closed in \( X \), and \( f \) is \( SC^* \)-\( gSC^* \)-closed, the image \( f(SC^*\text{-}\mathrm{cl}(J)) \) is \( gSC^* \)-closed in \( Y \).

Furthermore, by the properties of the \( SC^* \)-closure (cf. Lemma 1.10), we have:
\[
SC^*\text{-}\mathrm{cl}(f(J)) \subseteq SC^*\text{-}\mathrm{cl}(f(SC^*\text{-}\mathrm{cl}(J))) \subseteq N.
\]
Therefore, \( f(J) \) lies within a \( gSC^* \)-closed subset of \( N \), implies that \( f(J) \) is \( gSC^* \)-closed in \( Y \).

\subsection{Definition.} A function \(X \) to \(Y\) is termed \textbf{\( SC^* \)-irresolute} if $\forall$ \( SC^* \)-open set \( N \) in \( Y \), the preimage \( f^{-1}(N) \) is an \( SC^* \)-open set in \( X \).

\subsection{Proposition} 
Let \(X \) to \(Y\) be an open, \( SC^* \)-irresolute bijection. If \( I \subseteq Y \) is \( gSC^* \)-closed, then the preimage \( f^{-1}(I) \) is \( gSC^* \)-closed in \( X \).

\textbf{Proof}: Suppose \( M \) is an open set in \( X \) that contains \( f^{-1}(I) \). This implies \( I \subset f(M) \), and since \( f(M) \) is open in \( Y \), we have \( SC^* \)-\( cl(I) \subset f(M) \). Consequently, \( f^{-1}(SC^* \)-\( cl(I)) \subset M \). Given that \( f \) is \( SC^* \)-irresolute, \( f^{-1}(SC^* \)-\( cl(I)) \) remains \( SC^* \)-closed in \( X \). By applying \textbf{Lemma 1.9 $(i)$ $\&$ $(v)$}, it follows that \( SC^* \)-\( cl(f^{-1}(I)) \subset f^{-1}(SC^* \)-\( cl(I)) \subset M \). Hence, \( f^{-1}(I) \) is \( gSC^* \)-closed in \( X \).

\subsection{Theorem}
Let \( f: X \to Y \) $\&$ \( g: Y \to Z \) be functions. Then:

\begin{flushleft}
    \textbf{$(i)$} If \( f \) is a cont., surjection $\&$ the composition \( g \circ f: X \to Z \) is \( gSC^* \)-closed, then \( g \) is \( gSC^* \)-closed.\\
    \textbf{$(ii)$} If \( f: X \to Y \) and \( g: Y \to Z \) s.t. \( g \circ f \) is \( gSC^* \)-closed, $\&$ \( g \) both cont., and \( SC^* \)-\( gSC^* \)-closed, then \( g \circ f \) is \( gSC^* \)-closed.\\
    \textbf{$(iii)$} If \( g \circ f: X \to Y \) is closed and \( g \) is \( gSC^* \)-closed, then the composition \( g \circ f: X \to Z \) is \( gSC^* \)-closed.
\end{flushleft}

\textbf{Proof.}
\textbf{$(i)$} Suppose \( F \subseteq Y \) be closed set. If \( f \) is cont., the preimage \( f^{-1}(F) \) is closed in \( X \). Given that \( g \circ f \) is \( gSC^* \)-closed, the image \( g(f^{-1}(F)) = g(F) \) is \( gSC^* \)-closed in \( Z \). Thus, \( g \) preserves \( gSC^* \)-closedness and is \( gSC^* \)-closed.\\
\textbf{$(ii)$} This result follows directly from \textbf{Proposition 4.8}, which guarantees the preservation of \( gSC^* \)-closedness under the given conditions.\\
\textbf{$(iii)$} If \( g \circ f \) is a closed map and \( g \) is \( gSC^* \)-closed, then by definition of composition and preservation under closure, \( g \circ f \) is \( gSC^* \)-closed.

\subsection{Theorem.}  
Following conditions are equivalent in a topological space \( X \):  
\begin{enumerate}[(i)]
    \item \( X \) is \( SC^* \)-regular.
    \item Given any closed set \( F \in X\), a point \( x \notin F \), $\exists$ an \( SC^* \)-open set \( M \) $\&$ \( gSC^* \)-open set \( N \) s.t. \( x \in M \), \( F \subset N \), \( M \cap N = \emptyset \).
    \item Every subset \( J \) of \( X \), any closed set \( F \) with \( J \cap F = \emptyset \), $\exists$ an \( SC^* \)-open set \( M \) and a \( gSC^* \)-open set \( N \) satisfying \( J \cap M \neq \emptyset \), \( F \subset N \), $\&$ \( M \cap N = \emptyset \).
    \item For each closed set \( F \) in \( X \), the equality  
          \[
          F = \bigcap \{ SC^* \text{-} cl(N) \mid F \subset N, \, N \text{ is } SC^* \text{-open} \}
          \]
          holds.
\end{enumerate}

\textbf{Proof}:  

\textbf{$(i)$ \(\Rightarrow\) $(ii)$:}  
This follows directly since an \( SC^* \)-regular space ensures the separation of a point and a closed set by disjoint \( SC^* \)-open and \( gSC^* \)-open sets.  

\textbf{$(ii)$ \(\Rightarrow\) $(iii)$:}  
Let \( J \subset X\) and \( F \) a closed set s.t. \( J \cap F = \emptyset \). For any \( x \in J \), we have \( x \notin F \). By assumption, $\exists$ an \( SC^* \)-open set \( M \) and a \( gSC^* \)-open set \( N \) s.t. \( x \in M \), \( F \subset N \), \( M \cap N = \emptyset \). Since \( x \in J \), \( x \in M \), we conclude that \( J \cap M \neq \emptyset \), proving $(iii)$.  

\textbf{$(iii)$ \(\Rightarrow\) $(i)$:}  
Suppose \( F \) is a closed set in \( X \) and let \( x \notin F \). Since \( \{x\} \cap F = \emptyset \), we can apply condition $(iii)$, to obtain an \( SC^* \)-open set \( M \), a \( gSC^* \)-open set \( W \) s.t. \( x \in M \), \( F \subset W \), \( M \cap W = \emptyset \). Defining \( N = SC^* \text{-int}(W) \), it follows from a previous result that \( F \subset N \), \( N \) is \( SC^* \)-open, and \( M \cap N = \emptyset \). Hence, \( X \) is \( SC^* \)-regular.  

\textbf{$(i)$ \(\Rightarrow\) $(iv)$:}  
For any closed set \( F \) in \( X \), by a known theorem,  
\[
F \subset \bigcap \{ SC^* \text{-cl}(N) \mid F \subset N, \, N \text{ is } gSC^* \text{-open} \}.
\]  
Since every \( gSC^* \)-open set is \( SC^* \)-open, we obtain  
\[
F = \bigcap \{ SC^* \text{-cl}(N) \mid F \subset N, \, N \text{ is } SC^* \text{-open} \}.
\]  

\textbf{$(iv)$ \(\Rightarrow\) $(i)$:}  
Let \( F \) be a closed set in \( X \), \( x \notin F \). By $(iv)$, $\exists$ a \( gSC^* \)-open set \( W \) containing \( F \) s.t. \( x \notin SC^* \text{-cl}(W) \). Since \( F \) is closed, it follows from a \textbf{Lemma 1.10}, that \( F \subset SC^* \text{-int}(W) \). Setting \( N = SC^* \text{-int}(W) \), we get \( F \subset N \) and \( N \) is \( SC^* \)-open. Since \( x \notin SC^* \text{-cl}(W) \), we also have \( x \notin SC^* \text{-cl}(N) \). Defining \( M = X \setminus SC^* \text{-cl}(N) \), we obtain \( x \in M \), \( M \) is \( SC^* \)-open, and \( M \cap N = \emptyset \), proving that \( X \) is \( SC^* \)-regular.

\subsection{Definition} A function $f : X\rightarrow Y$ is said to be $SC^*$-open\cite{Chandrakala2022} if for each open set $M$ of $X$, $f(M)\in SC^*$-$O(Y)$.

\subsection{Theorem} Let \( f: X \to Y \) be a continuous function that is both \( SC^* \)-open and \( gSC^* \)-closed, and suppose \( f \) is surjective. Thus \( X \) is regular space, then \( Y \) is \( SC^* \)-regular.  

\textbf{Proof}:  
Consider any point \( y\) in \(Y\), an open set \( N \) in \( Y \) s.t. \( y \in N \). Since \( f \) is onto, $\exists$ \( x \in X \) s.t. \( f(x) = y \). Given that \( X \) is regular, we can find an open set \( M \) in \( X \) that satisfies  
\[
x \in M \subset \text{cl}(M) \subset f^{-1}(N).
\]  
Applying the function \( f \) to these inclusions, we obtain  
\[
y \in f(M) \subset f(\text{cl}(M)) \subset N.
\]  
Since \( f \) is \( SC^* \)-open, \( f(M) \) belongs to \( SC^*O(Y) \), $\&$ \( f \) is also \( gSC^* \)-closed, \( f(\text{cl}(M)) \) is \( gSC^* \)-closed in \( Y \). Consequently, 
\[
y \in f(M) \subset SC^*\text{-cl}(f(M)) \subset SC^*\text{-cl}(f(\text{cl}(M))) \subset N.
\]  
By an established result (\textbf{Theorem 4.12}), this confirms that \( Y \) satisfies the conditions of an \( SC^* \)-regular space.

\subsection{Definition} A function $f : X\rightarrow Y$ is said to be pre$SC^*$-open\cite{tomarsome2024} if for each $SC^*$-open set $M$ of $X$, $f(M)\in SC^*$-$O(Y)$.

\subsection{Theorem} 
Let \(X \) to \(Y\) be a cont., pre-\( SC^* \)-open, and \( SC^* \)-\( gSC^* \)-closed surjective function. If \( X \) an \( SC^* \)-regular space, then \( Y \) also \( SC^* \)-regular.  

\textbf{Proof}:  
Suppose that \( F \subset Y\), let \( y \) a point in \( Y \) but not in \( F \). The closed preimage \( f^{-1}(F) \subset X\). If \( X \) is \( SC^* \)-regular, so for any \( x \in f^{-1}(y) \), $\exists$ \( SC^* \)-open sets \( M \) and \( N \) in \( X \) s.t.:  
\[
x \in M, \quad f^{-1}(F) \subset N, \quad \text{and} \quad M \cap N = \emptyset.
\]  
Because \( F \) closed in \( Y \), $\exists$ an \( SC^* \)-open set \( W \) in \( Y \) s.t.:  
\[
F \subset W \quad \text{and} \quad f^{-1}(W) \subset N.
\]  
If \( f \) is pre-\( SC^* \)-open, so \( y = f(x) \in f(M) \),  \( f(M) \) is an \( SC^* \)-open set in \( Y \). Furthermore, from \( M \cap N = \emptyset \), it follows that  
\[
f^{-1}(W) \cap M = \emptyset \quad \Rightarrow \quad W \cap f(M) = \emptyset.
\]  
Ensures \( Y \) satisfies the separation property required for \( SC^* \)-regularity. Therefore, \( Y \) is \( SC^* \)-regular.

\section{Preservation theorems and other characterizations of $SC^*$-regular spaces}
\subsection{Theorem.} If \(X\) to \(Y\) is a cont., quasi-\( SC^* \)-closed, \( gSC^* \)-closed surj., function and \( X \) is \( SC^* \)-regular, \( Y \) is also $r$-space.  

\textbf{Proof}:  
Consider two disjoint closed sets \( G_1 \) and \( G_2 \) in \( Y \). If \( f \) is cont., their preimages, \( f^{-1}(G_1) \), \( f^{-1}(G_2) \), disjoint closed subsets of \( X \).  

If \( X \) is \( SC^* \)-regular, $\exists$ disjoint \( SC^* \)-open sets \( M_1 \), \( M_2 \) in \( X \) s.t.:  
\[
f^{-1}(G_1) \subset M_1, \quad f^{-1}(G_2) \subset M_2.
\]  
Define the sets in \( Y \) as follows:
\[
N_i = Y - f(X - M_i), \quad \text{for } i = 1,2.
\]  
Since \( f \) is a quasi-\( SC^* \)-closed and \( gSC^* \)-closed function, each \( N_i \) is an open set in \( Y \). Moreover, these sets satisfy:  
\[
G_i \subset N_i \quad \text{and} \quad f^{-1}(N_i) \subset M_i, \quad \text{for } i = 1,2.
\]  
As \( M_1 \) and \( M_2 \) are disjoint, their images under \( f \) also remain disjoint, ensuring:
\[
N_1 \cap N_2 = \emptyset.
\]  
This establishes that \( Y \) satisfies the separation condition for regularity, proving that \( Y \) is a regular space.

\subsection{Lemma.} 
If a subset \( J \) of a topological space \( X \) is \( gSC^* \)-open iff, $\forall$ closed set \( F \) contained in \( J \), we have \( F \subset SC^* \)-\(\operatorname{int}(J) \).

 \subsection{Theorem.} If \( X\) to \(Y\) be a closed and \( SC^* \)-\( gSC^* \)-cont., inject., function. If \( Y \) is \( SC^* \)-regular, \( X \) is also \( SC^* \)-regular.  

\textbf{Proof}:  
Consider two disjoint closed subsets \( V_1 \) and \( V_2 \) of \( X \). Since \( f \) is an injective function that is also closed, the images \( f(V_1) \) and \( f(V_2) \) are closed and disjoint in \( Y \).  

Given that \( Y \) is \( SC^* \)-regular, we can find disjoint \( SC^* \)-open sets \( N_1, N_2 \) in \( Y \),   
\[
f(V_i) \subset N_i, \quad \text{for } i = 1,2.
\]  
Because \( f \) is \( SC^* \)-\( gSC^* \)-cont., the preimages \( f^{-1}(N_1) \), \( f^{-1}(N_2) \) form \( gSC^* \)-open sets in \( X \), ensuring:  
\[
V_i \subset f^{-1}(N_i), \quad \text{for } i = 1,2.
\]  
So, define:  
\[
M_i = SC^*\text{-int}(f^{-1}(N_i)), \quad \text{for } i = 1,2.
\]  
Since \( SC^*\text{-int}(f^{-1}(N_i)) \) is an \( SC^* \)-open set in \( X \), we obtain:  
\[
M_i \in SC^*O(X), \quad V_i \subset M_i, \quad M_1 \cap M_2 = \emptyset.
\]  
Thus, \( X \) satisfies the separation condition required for \( SC^* \)-regularity, proving that \( X \) is \( SC^* \)-regular.   

\subsection{Corollary.}  
Let \( f: X \to Y \) be an injective function that is both closed and \( SC^* \)-irresolute. If \( Y \) is an \( SC^* \)-regular space, then \( X \) must also be \( SC^* \)-regular.

\textbf{Proof}: The result follows directly from the fact that every \( SC^* \)-irresolute function is also \( SC^* \)-\( gSC^* \)-continuous.

\subsection{Lemma.} A function \(X \) to \(Y\) is almost-\( gSC^* \)-closed iff $\forall$ subset \( I \) of \( Y \) $\&$ \( M \in RO(X) \) with \( f^{-1}(I) \subset U \), $\exists$ a \( gSC^* \)-open set \( N \) in \( Y \) s.t. \( I \subset N \) and \( f^{-1}(N) \subset M \).

\subsection{Lemma.} If \(X\) to \(Y\) is almost-\( gSC^* \)-closed, for any closed set \( G \) in \( Y \) any  \( M \in RO(X) \) containing \( f^{-1}(G) \), $\exists$ set \( N \in SC^*O(Y) \) satisfying \( G \subset N \) and \( f^{-1}(N) \subset M \).

\subsection{Theorem.}  
Let \( f: X \to Y \) be a surjective function that is continuous and almost \( gSC^* \)-closed. If \( X \) is a regular topological space, then \( Y \) is \( SC^* \)-regular.  

\textbf{Proof}:  
Consider two disjoint closed subsets \( M_1 \) and \( M_2 \) in \( Y \). If \( f \) is cont., their preimages \( f^{-1}(M_1) \), \( f^{-1}(M_2) \) closed sets in \( X \).  

If \( X \) is regular, then disjoint open sets \( U_1, U_2 \) in \( X \),  
\[
f^{-1}(M_i) \subset U_i, \quad \text{for } i = 1,2.
\]  
Define:  
\[
G_i = \operatorname{int}(\operatorname{cl}(U_i)), \quad \text{for } i=1,2.
\]  
Since regular open sets remain disjoint, we have:  
\[
G_1 \cap G_2 = \emptyset.
\]  

By \textbf{Lemma 5.6}, $\exists$ sets \( N_1, N_2 \in SC^*O(Y) \) s.t.:  
\[
M_i \subset N_i \quad \text{and} \quad f^{-1}(N_i) \subset G_i, \quad \text{for } i=1,2.
\]  
If \( f \) is a surjection and \( G_1 \cap G_2 = \emptyset \), it follows that:  
\[
N_1 \cap N_2 = \emptyset.
\]  
This satisfies the separation condition for \( SC^* \)-regularity in \( Y \), proving \( Y \) is \( SC^* \)-regular.

\section*{Conclusion}  
In this study introduced the concept of \( SC^* \)-regular spaces through the framework of \( SC^* \)-open sets. We investigated their connections with various generalizations of regularity, including regular, strongly, almost, softly, weakly, alpha, and generalized-regular spaces. Key characterizations and properties of \( SC^* \)-regular spaces were established, alongside an analysis of function types such as \( gSC^* \)-closed and \( SC^*g \)-closed mappings. Furthermore, we proved several preservation theorems that demonstrate how \( SC^* \)-regularity behaves under continuous, surjective, and specialized mappings. The results presented offer a foundation for further exploration in generalized topological structures.


\bibliographystyle{amsplain}
\bibliography{references}

\end{document}